\documentclass[12pt]{article}
\usepackage{widetext}
\usepackage{url}
\usepackage{quasilma}
\usepackage{hang}

\begin{document}

\def\G{\Gamma}
\def\f{\varphi}
\def\e{\varepsilon}
\def\TT{\mathop{\cal T}\nolimits}
\def\FF{\mathop{\cal F}\nolimits}
\def\GG{\mathop{\cal G}\nolimits}
\def\eql    {\mathop{=}      \limits}
\def\ul     {\mathop{\rm u}  \limits}
\def\ulo{\ul^{\scriptscriptstyle \circ}}
\def\suml   {\mathop{\sum}   \limits}
\def\bigcupl{\mathop{\bigcup}\limits}
\def\setminus#1#2{#1\!\smallsetminus\!#2}
\def\1n{1,\hdots,n}
\def\_#1{\mathop{\hspace{-2pt}^{}_{#1}}}
\newcommand{\x}[1]{}  
\def\hang{\indent}

\title{Matrix-Forest Theorems}

\author{PAVEL CHEBOTAREV\footnote{Corresponding author.
E-mail: {\tt pavel4e@gmail.com\x{, pchv@rambler.ru}}} {\small AND} ELENA SHAMIS}
\institution{Institute of Control Sciences of the Russian Academy of Sciences\\
65 Profsoyuznaya Street, Moscow 117997, Russian Federation}

\abstract{
\noindent The Laplacian matrix of a graph $G$ is $L(G)=D(G)-A(G)$,
where $A(G)$ is the adjacency matrix and $D(G)$ is the diagonal matrix
of vertex degrees. According to the Matrix-Tree Theorem, the number of
spanning trees in $G$ is equal to any cofactor of an entry of $L(G)$.
A rooted forest is a union of disjoint rooted trees. We consider the
matrix $W(G)=I+L(G)$ and prove that the $(i,j)$-cofactor of $W(G)$ is
equal to the number of spanning rooted forests of $G$, in which the
vertices $i$ and $j$ belong to the same tree rooted at $i$. The
determinant of $W(G)$ equals the total number of spanning rooted
forests, therefore the $(i,j)$-entry of the matrix $W^{-1}(G)$ can be
considered as a measure of relative ``forest-accessibility'' of the
vertex $i$ from $j$ (or $j$ from $i$). These results follow from
somewhat more general theorems we prove, which concern weighted
multigraphs. The analogous theorems for (multi)digraphs are
established. These results provide a graph-theoretic
interpretation to the adjugate of the Laplacian characteristic matrix.
}

\section{NOTE OF NOVEMBER 2023}

\noindent The first version of this arXiv preprint was the 31.03.1997 version of an unpublished manuscript of 12.04.1994, while the latter contained the same proofs of the same theorems, but only 18 bibliographic references. The differences of the present version are the exact\x{ bibliographic} reference \cite{ChShAT} (with the volume and pages), mentioning the version of \cite{ChShAtl} that is still available online, and this note.

\section{1. INTRODUCTION}

\noindent Let $G$ be a labeled graph on $n$ vertices with adjacency matrix
$A(G)=\left(a_{ij}\right)$. The Laplacian (the Kirchhoff or the
admittance) matrix of $G$ is the $n$-by-$n$ matrix
$L(G)=\left(\ell_{ij}\right)$ with $\ell_{ij}=-a_{ij}\; (j\neq i,\;\:
i,j=\1n)$ and $\ell_{ii}=\suml_{j\neq i}a_{ij}=-\suml_{j\neq
i}\ell_{ij}\; (i=\1n)$. According to the Matrix-Tree Theorem
attributed to Kirchhoff (for its history, see \cite{15}), any cofactor
of an entry of $L(G)$ is equal to the number of spanning trees of $G$.
Tutte (see \cite{19}) has generalized this theorem to weighted
multigraphs and multidigraphs.  Bapat and Constantine \cite{1}
presented a version for graphs in which each edge is assigned a color.
Merris \cite{13} proposed an ``edge version'' of the Matrix-Tree
Theorem and Moon \cite{MoonLAMA} generalized it. Forman \cite{5}
considered the Kirchhoff theorem in a more general context of vector
fields.

Another trend of literature studies the characteristic polynomial and
the spectrum of the Laplacian matrix. For review of this literature we
refer to \cite{6,7,14}. We would like to mention here the
research by Kelmans, who had published in 1965--1967 a series of
results on the Laplacian characteristic polynomial and spectrum (see
\cite{9,10}, and the references therein), some of which were
rediscovered later by other writers.

In \cite{10} Kelmans and Chelnokov have shown that the coefficients of
the Laplacian characteristic polynomial can be represented through the
numbers of spanning forests of $G$ with fixed numbers of components.
This result is closely related to the theorems in this paper and we use
it in our proofs. Another close result has been obtained by Fiedler and
Sedl\'{a}\v{c}ek \cite{4} (see Lemma~3 in the Appendix) and generalized
in \cite{Chen,Cha,MoonDM}.

We consider the matrix $W(G)=I+L(G)=-Z(-1,G)$, where $Z(\lambda,G)=
\lambda I-L(G)$ is the Laplacian characteristic matrix of $G$ and $I$
is the identity matrix. It turns out that $W(G)$ can be used for
counting spanning rooted forests of $G$ (a rooted forest is a union of
disjoint rooted trees) somewhat similarly to as $L(G)$ can be used to
count spanning trees. Namely, the determinant of $W(G)$ is equal to the
number of all spanning rooted forests of $G$, and the cofactor
$W^{ij}(G)$ is equal to the number of those spanning rooted forests,
such that $i$ and $j$ belong to the same tree rooted in $i$. This is a
simple consequence of Theorems~5 and 6 in Section 3.

Theorems~3 and 4 consider an arbitrary multidigraphs $\G$ and provide
an analogous relation between the Kirchhoff matrix of $\G$ and the
numbers of spanning {\it diverging} forests of $\G$. These results
allow us to consider the matrices $W^{-1}(G)$ and $W^{-1}(\G)$ as
matrices of {\it relative forest-accessibilities} in the multigraph $G$
and the multidigraph $\G$.

It can be interesting to compare Theorems~3--7 with the
corresponding results on the adjacency characteristic matrix (see
\cite[Subsections 1.9.1 and 1.9.5]{2} or the original articles by
Kasteleyn and Ponstein cited therein, and also \cite{17}). About
counting forests see \cite{3,8,15}. Liu and Chow \cite{11} obtained a
rather complicated expression for the number of $k$-component spanning
forests of a graph through the principal minors of the Laplacian
matrix.  Myrvold \cite{16} gave a simpler graph-theoretic proof of some
version of their result and discovered a polynomial algorithm for
calculating this number of $k$-component spanning forests. The ideas of
her proof are similar to those of Kelmans and Chelnokov \cite{10}.

In the following section, we give some graph-theoretic notation and
statements of the Matrix-Tree Theorem for weighted multigraphs and
multidigraphs, which will help us to formulate and prove our results.

\section{2. PRELIMINARIES}

\noindent Let us remind some necessary graph-theoretic notions. A {\it subgraph}
of a multigraph $G$ is a multigraph all of whose vertices and edges
belong to the vertex and edge sets of $G$. A {\it spanning subgraph} of
$G$ is a subgraph of $G$ whose vertex set coincides with the vertex set
of $G$. A {\it forest} is a cycleless graph. A {\it tree} is a
connected forest. A {\it rooted tree} is a tree with one marked vertex
called a {\it root}. Formally, the rooted tree is a pair $(T,r)$, where
$T$ is the tree and $r$ is its vertex. A {\it component} of a
multigraph $G$ is any maximal (by inclusion) connected subgraph of $G$.
A {\it rooted forest} can be defined as a forest with one marked vertex
in each component. Obviously, a rooted forest is a union of disjoint
rooted trees.

A {\it path} from vertex $i$ to vertex $j$ in a multidigraph $\G$ is an
alternating sequence of distinct vertices and arcs having each arc
directed from the previous vertex to the next one; $i$ is the first
vertex, $j$ is the last one. Sometimes we will consider a path as a
subgraph of $\G$. A digraph is called a {\it directed tree} (a {\it
directed forest}) if the graph obtained from it by replacement of all
its arcs by edges is a tree (a forest). The definitions for {\it
directed rooted tree} and {\it directed rooted forest} are analogous to
the definitions of rooted tree and rooted forest (we will omit the word
``directed'' while talking about subgraphs of $\G$). A {\it diverging
tree} is a directed rooted tree, containing paths from the root to all
other vertices. A {\it diverging forest} is a directed rooted forest,
all whose components are diverging trees.

The Matrix-Tree Theorem asserts that for any graph $G$, all cofactors
of entries of $L(G)$ are equal to each other and their common value is
the number of spanning trees in $G$.

Now suppose $G$ is a weighted multigraph with vertex set
$V(G)=\{\1n\}$ and suppose $\e_{ij}^m$ is the weight of the
$m$th edge between vertices $i$ and $j$ in $G$. This weight will be
also referred to as a {\it conductance} of the edge. However, we will
not forbid $\e_{ij}^m$ to be negative. The Kirchhoff matrix of $G$ is
the $n$-by-$n$ matrix $L=L(G)=\left(\ell_{ij}\right)$ with
$\ell_{ij}=-\suml_{m=1}^{a\_{ij}}\e_{ij}^m\;\: (j\neq i,\;\:
i,j=\1n)$ and $\ell_{ii}=-\suml_{j\neq i}\ell_{ij}\;\:
(i=\1n)$, where $a_{ij}$ is the number of edges between $i$ and
$j$. Denote by $L^{ij}$ the cofactor of $\ell_{ij}$ in $L$. The product
of the conductances of all edges belonging to a subgraph $H$ of the
multigraph $G$ will be referred to as the {\it weight} or
{\it transmission coefficient} of $H$ and denoted by $\e(H)$.
The weight of a subgraph without edges is assumed to be 1. For every
nonempty set of subgraphs $\GG$, its weight is defined as follows:
$\e(\GG)=\suml_{H\in\GG}\e(H).$ Set the weight of the empty set to be
zero. Let $\TT(G)=\TT$ be the set of all spanning trees of $G$.

Tutte's \cite{19} generalization of the Matrix-Tree Theorem can be
formulated as follows.

\statement{Theorem~1}{(Matrix-Tree Theorem for weighted
multigraphs)}{For any weighted multigraph $G$ and for any $i,j\in
V(G)$,  $L^{ij}=\e(\TT)$.}

Tutte \cite{19} also developed a parallel theory for multidigraphs.

Let $\G$ be a multidigraph with vertex set $V(\G)=\{\1n\}$ and
suppose $\e_{ij}^m$ is the weight (or the conductance) of the $m$th arc
from $i$ to $j$ in $\G$. The Kirchhoff matrix $L(\G)$ of $\G$ is the
$n$-by-$n$ matrix $L=L(\G)=\left(\ell_{ij}\right)$ with
$\ell_{ij}=-\suml_{m=1}^{a\_{ji}}\e_{ji}^m\:\: (j\neq i,\;\:
i,j=\1n)$ and $\ell_{ii}=-\suml_{j\neq i}\ell_{ij}\:\:
(i=\1n)$, where $a_{ji}$ is the number of arcs from $j$ to $i$
in $\G$. Notice that $\ell_{ii}$ is the total conductance of the arcs
{\it converging} to $i$. The definition for weight of
a subgraph of $\G$ is analogous to the corresponding definition for
multigraphs. Suppose $\TT^i$ is the set of spanning trees of $\G$
diverging from $i$.

\statement{Theorem~2}{(Matrix-Tree Theorem for weighted
multidigraphs)}{For any weighted multidigraph $\G$ and for any $i,j\in
V(\G)$, $L^{ij}=\e(\TT^i)$.}

Note that in the directed case entries in different rows of $L$
may have different cofactors, but all the entries of one row have
equal cofactors.

For simplicity, Tutte formulates this theorem (as well as the previous
one) only for diagonal cofactors $L^{ii}$. The ``directed'' Matrix-Tree
Theorem concerning arbitrary $L^{ij}$ is given in Harary and Palmer
\cite{8}.

In the following section, we give somewhat analogous theorems on
spanning converging forests of a multidigraph $\G$ and on spanning
rooted forests of a multigraph $G$.

\section{3. MATRIX-FOREST THEOREMS}

\noindent Consider the matrices $W(\G)=I+L(\G)$ and $W(G)=I+L(G)$. $W^{ij}(\G)$
and $W^{ij}(G)$ will denote the cofactors of the $(i,j)$-entries of
$W(\G)$ and $W(G)$.

Suppose $\FF(\G)=\FF$ is the set of all spanning diverging forests of
$\G$ and $\FF^{i\rightarrow j}(\G)= \FF^{i\rightarrow j}$ is the set of
those spanning diverging forests of $\G$, such that $i$ and $j$ belong
to the same tree diverging from $i$. Let $W=W(\G),\;W^{ij}=W^{ij}(\G)$.

\statement{Theorem~3}{}{For any weighted multidigraph $\G,\;\: \det W=
\e(\FF)$.}

\secondstatement{Theorem~4}{}{For any weighted multidigraph $\G$ and
for any $i,j\in V(\G),\;\: W^{ij}=\e(\FF^{i\rightarrow j})$.}

As usual, these theorems have dual counterparts concerning {\it
converging} forests. Theorems~3 and~4 can be derived in the shortest
way from one version of Chaiken's result \cite{Cha}, namely, by
putting $U=W=\varnothing$ and then $U=\{i\},\: W=\{j\}$ in the first
formula in page~328 (cf.\ \cite[Theorem~3.1]{MoonDM}). In the Appendix
of this paper, we give another proof which demonstrates some
interesting relations of Matrix-Forest Theorems to the results in
\cite{4,9,10,12}.

Suppose $\FF(G)=\FF$ is the set of all spanning rooted forests of a
weighted multigraph $G$ and $\FF^{ij}(G)= \FF^{ij}$ is the set of those
spanning rooted forests of $G$, such that $i$ and $j$ belong to the
same tree rooted in $i$. Let $W=W(G),\;W^{ij}=W^{ij}(G)$.

\statement{Theorem~5}{}{For any weighted multigraph $G,\;\: \det W=
\e(\FF)$.}

\secondstatement{Theorem~6}{}{For any weighted multigraph $G$ and for
any $i,j\in V(G),\;$ $W^{ij}=\e(\FF^{ij})$.}

As the matrix $W$ of a weighted multigraph is symmetrical, Theorem~6
remains true if we replace $\FF^{ij}$ by $\FF^{ji}$ in the right-hand
side. In the Appendix, Theorems~5 and~6 are derived from Theorems~3
and~4.

If the weights $\e_{ij}^m$ are non-negative, then by Theorems~3 and 5,
the matrix $W$ is non-singular. If the matrix $W^{-1}$ exists, we will
denote it by $Q=\left(q\_{ij}\right)$ (both for a weighted
multidigraph $\G$ and for a weighted multigraph $G$). Then $Q=(\det
W)^{-1}W^*$, where $W^*=\left(W^{ij}\right)^{
\scriptscriptstyle \intercal
}$ is the adjugate
of $W$.  Theorems~3-6 imply the following main theorem.

\statement{Theorem~7}{}{1. For any weighted multidigraph $\G$, if the
matrix $Q=W^{-1}$ exists, then \newline
$q\_{ij}=\e(\FF^{j\rightarrow i})\Big/\e(\FF),\;\;i,j=\1n$.\newline
2. For any weighted multigraph $G$ if the matrix $Q=W^{-1}$ exists,
then \newline
$q\_{ij}=\e(\FF^{ji})\Big/\e(\FF),\;\;i,j=\1n$.}

It can be seen that $\suml_{j=1}^n q\_{ij}=1\;\; (i=\1n)$ both
for directed and undirected weighted multigraphs. This follows, for
example, from the facts that for any $i\in V(\G),$ the sets
$\FF^{j\rightarrow i}\;(j=\1n)$ are non-overlapping and
$\bigcupl_{j=1}^n\FF^{j\rightarrow i}=\FF$ (respectively, for any $i\in
V(G),$ the sets $\FF^{ji}\;(j=\1n)$ are non-overlapping and
$\bigcupl_{j=1}^n\FF^{ji}=\FF$).

If the weights of all arcs of $\G$ (of all edges of $G$) are ones,
Theorems~3--7 tell us about the {\it numbers} of corresponding spanning
forests (which are equal to their summary TC's in this case).

Theorem~7 allows us to consider the matrix $Q=W^{-1}$ as the matrix
of {\it relative forest-accessibilities} of the vertices of $\G$ (or
$G$).

Theorems~3 and~4 were formulated in \cite{ChShAtl} and Theorems~5
and~6 in \cite{ShaAMS}. The latter results were used in
\cite{18} for the analysis of one method of preference aggregation.
That paper implicitly contains proofs of these theorems (in the
case of equal weights $\e_{ij}^m$), different from the proofs given
here. In \cite{ChShAT} we analyze the properties of relative
forest-accessibilities and exploit them to introduce a new family of
sociometric indices.

Theorems~5--7 were used in \cite{18} for the analysis of one method of
preference aggregation. That paper contains implicitly the proofs of
these theorems (in the case of equal weights $\e_{ij}^m$), different
from the proofs given here.

\remark{}Lemma~5 in the Appendix provides an interpretation for the
adjugate of the characteristic matrix of $-L(\G)$. Replacing $\e(F)$
by $(-1)^{d(F)}\e(F)$ in (\ref{2}) ($d(F)$ is the number of arcs in
$F$), we obtain a representation for the adjugate of the characteristic
matrix of $L(\G)$. To get analogous representations in the undirected
case, it suffices to replace $\FF_{\f\cup\{i\}}^{i\rightarrow j}$ by
$\FF_{\f\cup\{i\}}^{ij}$ in (\ref{2}) ($\FF_{\f\cup\{i\}}^{ij}=
\FF^{ij} \cap \FF_{\f\cup\{i\}}$, and $\FF_{\f\cup\{i\}}$ is the set of
spanning rooted forests of $G$, having $|\f\cup\{i\}|$ components
rooted in the vertices of $\f\cup\{i\}$). The latter is obvious by the
argument used in the proof of Theorems~5 and~6.

\APPENDIX

Prior to proving Theorems~3--6 we introduce some additional notation
and prove several lemmas.

\hang $p(\lambda)$ ~is the characteristic polynomial of the matrix
$-L=-L(\G)$;

\hang $W_\lambda=\lambda I+L(\G)$ ($\lambda$ is a real number);

\hang $E=E(\Gamma)$ ~is the arc set of $\Gamma$;

\hang $\varphi$ ~is a subset of $V=V(\G)=\{\1n\}$;

\hang $L_{-\varphi}(\Gamma)=L_{-\varphi}$ ~is the matrix obtained from
$L(\Gamma)$ by deleting the rows and columns corresponding to the
vertices of $\varphi$; we will use the analogous expression $U_{-\psi}$
for an arbitrary $n$-by-$n$ matrix $U$ and $\psi\subseteq
\{\1n\}$.

\hang $\Gamma_{\varphi}$ ~is the weighted multidigraph obtained from
$\Gamma$ by identifying all the vertices of $\varphi$;

\hang $\varphi^*$ ~is the vertex of $\Gamma_{\varphi}$ being a result
of this identification; any arc incident to some vertex of $\varphi$ in
$\Gamma$ have the corresponding arc incident to $\varphi^*$ in
$\Gamma_{\varphi}$;

\hang ``$\G$-tree'' ~is a spanning diverging tree of $\G$;

\hang ``$\G$-forest'' ~is a spanning diverging forest of $\G$;

\hang $\TT_{\f^*}$ ~is the set of $\G_\f$-trees diverging from $\f^*$;
if $\f=\varnothing$, we set $\TT_{\f^*}=\varnothing$;

\hang $\FF_\f$ ~is the set of $\G$-forests with $|\f|$ components that
diverge from the vertices of $\f$;

\hang $\FF^{j \rightarrow i}_\f=\FF^{j \rightarrow i}\bigcap\FF_\f$ (if
$j\notin\f$, then $\FF^{j \rightarrow i}_\f=\varnothing$);

\statement{Lemma~1}{}{Let $\G_1$ and $\G_2$ be weighted multidigraphs
with the same vertex set. Suppose that the arc set of $\G_2$ can be
obtained from that of $\G_1$ by replacing some arc $(i,j)$ {\rm(}with
some weight $\e_{ij}${\rm)} by two arcs from $i$ to $j$ with the
weights $\e'_{ij}$ and $\e''_{ij}$ such that $\e'_{ij}+\e''_{ij}=
\e_{ij}$.
Then\newline
{\rm(}i{\rm)}~~$W(\G_1)=W(\G_2)${\rm;}\newline
{\rm(}ii{\rm)}~~for any vertices $\alpha$ and $\beta$, the value
$\e(\FF^{\alpha\rightarrow\beta})$ is the same for $\G_1$
and $\G_2$.}

\proof{}({\it i}) is obvious. ({\it ii}) holds since for any
$F\in\FF^{\alpha\rightarrow\beta}(\G_1)$ there are two corresponding
forests in $\FF^{\alpha\rightarrow\beta}(\G_2)$ with the same summary
weight.\endproof

\medskip
Based on Lemma~1, we conclude that it suffices to prove Theorems~3 and
4 only for weighted {\it digraphs}. Thus, we will assume that $\G$ has
no multiple arcs.

The following three lemmas are directed (and weighted) counterparts of
certain results of Kelmans \cite{9} and Kelmans and Chelnokov \cite{10}.

\statement{Lemma~2}{}{Let $\f\subseteq V$. Then, in terms of the
notation above, $\det L_{-\f}=\e(\TT_{\f^*}).$}

\proof{}If $\f=\varnothing$, we have zero in both sides of the
equality. For $\f\neq\varnothing,$ let $L(\G_\f)$ be the Kirchhoff
matrix of $\G_\f$. Suppose $L_{-\{\f^*\}}(\G_\f)=L_{-\f^*}(\G_\f)$ is
the matrix obtained from $L(\G_\f)$ by deleting the row and column
corresponding to $\f^*$.  Then the desired equality is valid since by
Theorem~2, $\det L_{-\f^*}(\G_\f)=\e(\TT_{\f^*})$, and
$L_{-\f^*}(\G_\f)=L_{-\f}$.\endproof

\statement{Lemma~3}{(Fiedler and Sedl\'{a}\v{c}ek \cite{4}, cf.\
\cite{1,5})}
{For any $\f\subseteq V$, $\det L_{-\f}=\e(\FF_\f).$}

We are proving Lemma~3 here, since this proof is very short.

\medskip
\proof{}By Lemma~2, it suffices to prove the equality
$\e(\FF_\f)=\e(\TT_{\f^*})$, which holds
for any $\f\neq\varnothing$ since identifying the vertices of $\f$
transforms any $\G$-forest belonging to $\FF_\f$ into a $\G_\f$-tree
diverging from $\f^*$, and this correspondence is one-to-one. If
$\f=\varnothing$, we have zero in both sides.\endproof

\medskip
Let $p(\lambda)=\det (\lambda I+L)=\suml^n_{k=0}c_k\lambda^k$ be
the characteristic polynomial of $-L$.

\statement{Lemma~4}{}{$c_k=\suml_{\textstyle{\f\subseteq V \atop
|\f|=k}}\e(\FF_\f),\;\;k=0,\hdots,n.$}

\proof{}In view of Lemma~3, this follows from the fact that $c_k$ is
equal to the sum of degree $k$ principal minors of $L$.\endproof

\bigskip
\proof{ of Theorem~3}Using Lemma~4, we have
$$
\det W=\det(I+L)=p(1)
=
\suml_{k=0}^n\suml_{\textstyle{\f\subseteq V \atop
|\f|=k}}\e(\FF_\f)
=\suml_{\f\subseteq V}\e(\FF_\f)
=\e(\FF).\endproofmath
$$

Suppose $W\_\lambda=\lambda I+L$ and
\begin{equation}
W_\lambda^{ij}=\suml_{k=0}^{n-1}b_k\lambda^k,\qquad i,j\in V
\label{1}
\end{equation}
is the cofactor of the $(i,j)$-entry of $W_\lambda$.

\statement{Lemma~5}{}{In terms of the notation above,
\begin{equation}
b_k=\suml_{\textstyle{\f\subseteq\setminus{V}{\{i,j\}} \atop
|\f|=k}}\e(\FF_{\f\cup\{i\}}^{i\rightarrow j}),\qquad
k=0,\hdots,n-1.
\label{2}
\end{equation}
}
\proof{}It is easy to see that
\begin{equation}
b_k=\suml_
{\textstyle{\f\subseteq\setminus{V}{\{i,j\}} \atop|\f|=k}}
L_{-\f}^{ij} \quad (k=0,\hdots,n-1),
\label{3}
\end{equation}
where $L_{-\f}^{ij}$ is the cofactor in the matrix $L_{-\f}$ of the
$(i,j)$-entry of $L$.

1. $i\neq j.$ To obtain an expression for $L_{-\f}^{ij}$, we will use a
theorem by Maybee (see \cite{12}), which can be formulated as follows.
For any $n$-by-$n$ matrix $U=(u_{ij})$, the representation of a
cofactor $U^{ij}$,
\begin{equation}
U^{ij}=\suml_k\e(P_k^{i\rightarrow j})\det U_{-\psi\_k},
\label{4}
\end{equation}
is valid for $i\neq j$, where $P_k^{i\rightarrow j}$ is the $k$th path
from $i$ to $j$ in an arbitrary weighted digraph $\G(U)$ (with vertex
set $\{\1n\}$) connected with $U$ by the following relations:
\newline
-- if $i\neq j$ and $u_{ij}\neq 0$, then the arc $(j,i)$ belongs to
arc set $E(\G(U))$ and has the weight $(-u_{ij})$;\newline
-- if $i\neq j$ and $u_{ij}=0$, then the arc $(j,i)$ has zero
weight or $(j,i)\notin E(\G(U))$.\newline
$\psi\_k$ in (\ref{4}) denotes the set of the vertices entering
$P^{i\rightarrow j}_k$.

Notice that matrix $L$ and weighted digraph $\G$ satisfy these
conditions (recall that by our assumption, $\G$ has no multiple arcs).
Therefore, these conditions are obeyed for $L_{-\f}$ and the subgraph
of $\G$ induced on the vertex subset $\setminus{V}{\f}$. Hence, by
(\ref{4}) and Lemma~3, we have
\begin{equation}
L_{-\f}^{ij}=\suml_k \e(P_k^{i\rightarrow j})
\det L_{-(\f\cup\psi\_k)}=\suml_k \e(P_k^{i\rightarrow j})
\e(\FF_{\f\cup\psi\_k}).
\label{5}
\end{equation}

Note that if $F\in\FF_{\f\cup\psi\_k}$ then the union of
$P_k^{i\rightarrow j}$ and $F$ belongs to $\FF_{\f\cup\{i\}}
^{i\rightarrow j}$. On the other hand, any $F'\in \FF_{\f\cup\{i\}}
^{i\rightarrow j}$ can be uniquely decomposed into a union of certain
$P_k^{i\rightarrow j}$ and $F\in\FF_{\f\cup\psi\_k}.$ Therefore,
(\ref{5}) implies
\begin{equation}
L_{-\f}^{ij}=\e(\FF_{\f\cup\{i\}}^{i\rightarrow j}).
\label{6}
\end{equation}

2. $i=j.$ Lemma~3 implies $L_{-\f}^{ij}=L_{-(\f\cup\{i\})}
=\e(\FF\_{\f\cup\{i\}}).$ Since
$\FF\_{\f\cup\{i\}}=\FF_{\f\cup\{i\}}^{i\rightarrow i}$, we have
(\ref{6}) as well.

Now (\ref{6}) and (\ref{3}) yield (\ref{2}).\endproof

\bigskip
\proof{ of Theorem~4}It suffices to put $\lambda=1$ in (\ref{1}) and use
Lemma~5:

$$
W^{ij}
=W_1^{ij}=\suml_{k=0}^{n-1}b_k
=\suml_{k=0}^{n-1}\suml_
{\textstyle{\f\subseteq\setminus{V}{\{i,j\}} \atop|\f|=k}}
\e(\FF_{\f\cup\{i\}}^{i\rightarrow j})
$$
$$
\phantom{1}=\suml_
{\f\subseteq\setminus{V}{\{i,j\}}}
\e(\FF_{\f\cup\{i\}}^{i\rightarrow j})
=\e(\FF^{i\rightarrow j}).\endproofmath
$$

\medskip
\proof{ of Theorems~5 and 6}Let $G$ be an arbitrary weighted {\it
graph} (by undirected counterpart of Lemma~1, we assume that $G$ has no
multiple edges). Replace every edge of $G$, having, say, a weight $\e$,
by two opposite arcs with the weight $\e$. The weighted digraph we
obtain has the same Kirchhoff matrix as $G$. The desired statements
follow from the fact that there exists a natural one-to-one
correspondence between rooted forests of $G$ and diverging forests of
$\G$.\endproof

\end{document}